\documentclass{amsart}

\newcommand{\bt}{\begin{Theorem}}
\newcommand{\et}{\end{Theorem}}
\newcommand{\bi}{\begin{itemize}}
\newcommand{\ei}{\end{itemize}}
\newcommand{\bea}{\begin{eqnarray}}
\newcommand{\ba}{\begin{array}}
\newcommand{\eea}{\end{eqnarray}}
\newcommand{\ea}{\end{array}}

\newcommand{\f}{\frac}
\newcommand{\s}{\sigma}

\newcommand{\Si}{\Sigma}
\newcommand{\iy}{\infty}
\newcommand{\what}{\widehat}
\newcommand{\lgra}{\longrightarrow}

\newtheorem{Definition}{Definition}[section]
\newtheorem{Theorem}[Definition]{Theorem}

\newtheorem{Corollary}[Definition]{Corollary}

\newtheorem{Remark}[Definition]{Remark}

\newcommand{\be}{\begin{equation}}
\newcommand{\ee}{\end{equation}}

\newcommand{\newsection}{\setcounter{equation}{0}}

\newcommand{\R}{\mathbb R}%
\newcommand{\C}{\mathbb C}%
\newcommand{\Z}{\mathbb Z}%

\begin{document}
\author[Sarkar]{Rudra P. Sarkar}
\address[R. P. Sarkar]{Stat-Math Unit, Indian Statistical
Institute, 203 B. T. Road, Calcutta 700108, India, E-mail:
rudra@isical.ac.in}
\author[Sengupta]{Jyoti Sengupta}
\address[J. Sengupta]{School of Mathematics, Tata Institute of Fundamental
Research, Homi Bhabha Rd.,  Mumbai 400005, India, E-mail:
sengupta@math.tifr.res.in}
\title[Beurling's Theorem and characterization of heat kernel]
{Beurling's Theorem and characterization of heat kernel for
Riemannian Symmetric spaces of noncompact type}
\begin{abstract}
We prove  Beurling's theorem for rank 1 Riemmanian symmetric
spaces and relate it to the characterization of the heat kernel of
the symmetric space.
\end{abstract}
\maketitle
\section{Introduction}
The uncertainty principle in harmonic analysis reflects the
inevitable trade-off between the function and its Fourier
transform as it says that both of them cannot decay very rapidly.
This principle has several quantitative versions which were proved
by Hardy, Morgan, Gelfand-Shilov, Cowling-Price etc. (see
\cite{HJ}, \cite{FS}, \cite{Th1} and the references there in). In
more recent times H\"ormander  (see \cite{Hor}) proved the
following theorem which is the strongest theorem in this genre in
the sense that it implies the theorems of Hardy, Morgan,
Gelfand-Shilov  and Cowling-Price.

\begin{Theorem}{\em (H\"ormander 1991)}
Let $f\in L^1(\R)$. Then
$$\int_\R\int_\R |f(x)||\what{f}(y)|e^{|x||y|}dxdy<\infty$$
implies $f= 0$ almost everywhere. \label{thm-Hormander}
\end{Theorem}

H\"ormander attributes this theorem to A. Beurling.

As is well-known in physics, the uncertainty in the momentum is
smallest, for a given uncertainty in the position, if the wave
function is the Gaussian $e^{-\frac {x^2}{4t}}$. In harmonic
analysis this means that the trade-off is optimal when the
function is Gaussian. The quantitative versions of the uncertainty
principle also accommodate this optimal situation. The above
theorem of H\"ormander was further generalized in \cite{BDJ} which
takes care of this aspect of uncertainty:

\begin{Theorem}{\em (Bonami, Demange, Jaming 2003)}
Let $f\in L^2(\R)$ and $N\ge 0$. Then
$$\int_\R\int_\R\frac{|f(x)||\what{f}(y)|}{(1+|x|+|y|)^N}e^{|x||y|}dxdy<\infty$$
implies $f(x)=P(x)e^{-tx^2}$ where $t>0$ and $P$ is a polynomial
with $\deg P<\frac{N-1}2$. \label{thm-BDJ}
\end{Theorem}
We will refer to theorem \ref{thm-BDJ} simply as {\em Beurling's
theorem} for the sake of brevity.

The aim of this article is to prove the analogue of theorem
\ref{thm-BDJ} for Riemannian symmetric spaces $X$ of the
noncompact type which have rank 1.   We recall that such a space
is of the form $G/K$ where $G$ is a noncompact connected
semisimple Lie group of real rank $1$ with finite centre and
$K\subset G$ is a maximal compact subgroup.

The precise statement of the theorem and its proof appear in
section 3. In section 4 we have showed that the estimate
considered in the main theorem is the sharpest possible. In
section 5 we have indicated how the theorems of Hardy, Morgan,
Gelfand-Shilov, Cowling-Price etc on symmetric spaces follow from
our Beurling's theorem. The mutual dependencies of these
uncertainty theorems can be schematically displayed as follows:
$$\ba{ccccc}\mbox{ Beurling's
}&\Rightarrow&\mbox{ Gelfand-Shilov }&\Rightarrow&\mbox{ Cowling-Price }\\
&&\Downarrow&&\Downarrow\\
&&\mbox{ Morgan's }&\Rightarrow&\mbox{ Hardy's }\ea$$ This shows
that Beurling's theorem is the {\em Master theorem}. Some of the
latter theorems (which follow from Beurling's) were proved
independently on symmetric spaces in recent years by many authors
(see \cite{SiSu, CSS, EEKK, NR1, NR2, Sengupta2, Th1, RS} etc.).

After completing this work we had the opportunity to see Demange's
thesis (\cite{Demange}) in which he further generalized theorem
1.2 (see theorem 6.1). In section 6 we have given the appropriate
analogue of Demange's theorem on symmetric spaces.

\section{Notation and Preliminaries}\newsection
The pair $(G,K)$ is as described in the introduction. We let
$G=KAN$ denote a fixed Iwasawa decomposition of $G$. Let
$\mathfrak g$, $\mathfrak k$, $\mathfrak a$ and  $\mathfrak n$
denote the Lie algebras of $G$, $K$, $A$ and $N$ respectively. We
recall that dimension of $\mathfrak a=1$. We choose and keep fixed
throughout a system of positive restricted roots, which we denote
by $\Sigma^+$. Let $\gamma\in \Sigma^+$ denote the unique simple
root and let $H_\gamma\in \mathfrak a$ be the dual basis of
$\mathfrak a$. Using $\gamma$ (respectively $H_\gamma$) we can
identify $\mathfrak a^*$ (respectively $\mathfrak a$) with $\R$.
The complexification $\mathfrak a^*_\C$ of $\mathfrak a^*$ can
then be identified with $\C$. Under this correspondence the
half-sum of the elements of $\Sigma^+$, denoted by $\rho$
corresponds to the real number $\frac 12 (m_\gamma+2m_{2\gamma})$
where $m_\gamma$ (respectively $m_{2\gamma}$) is the multiplicity
of the root $\gamma$ (respectively $2\gamma$). We will frequently
identify $\rho$ with this positive real number without further
comment. Furthermore, the positive Weyl chamber $\mathfrak
a_+\subset \mathfrak a$ (respectively $\mathfrak a_+^*\subset
\mathfrak a^*$) gets identified under this correspondence with the
set of positive real numbers. We let $\exp t H_\gamma=a_t\in A$
for $t\in \R$. This identifies $A$ with $\R$. Let $H:G\lgra
\mathfrak a$ be the Iwasawa projection associated to the Iwasawa
decomposition, $G= KAN$. Then $H$ is left $K$-invariant and right
$MN$-invariant where $M$ is the centraliser of $A$ in $K$. For
$\lambda\in \mathfrak a^*$ (respectively $H\in \mathfrak a$) we
denote by $\lambda^+$ (respectively $H^+$) the unique Weyl
translate of $\lambda$ (respectively $H$) that  belongs to the
closure of the positive Weyl chamber $\mathfrak a^*_+$
(respectively $\mathfrak a_+$). We have
$\lambda^+(H^+)=|\lambda(H)|$ where $|r|$ denotes the modulus of
the real number $r$.  Note that the Weyl group is isomorphic to
$\Z_2$. The unique nontrivial element of the Weyl group takes an
element $\lambda\in \mathfrak a^*\equiv \R$ (respectively $H\in
\mathfrak a$) to $-\lambda$ (respectively $-H$). Therefore
$\lambda^+$ (respectively $H^+$) corresponds to $|\lambda|$
(respectively $|H|$) under the above identification of $\mathfrak
a^*$ (respectively $\mathfrak a$) and $\R$.

We have the $\mathfrak a$-valued inner product of Helgason $A(x,
k)$ on $X\times K$ defined by $A(x, k)=-H(x^{-1}k), x\in X, k\in
K$. Note that $A$ descends to a function, also denoted by
$A:X\times K/M\lgra \mathfrak a$, since $H$ is right
$M$-invariant.

We fix a left $G$-invariant measure $dx$ on $X$ and a left
invariant Haar measure $dg$ on $G$ such that for a nice function
$f$ on $X$, $\int_Xf(x)dx=\int_Gf(g)dg$. Here in the right hand
side $f$ is considered as a right $K$-invariant function on $G$.
While dealing with functions on $X$, we may slur over the
difference between the two measures and denote  both by $dx$. Let
$dn$ (respectively $da$) be a fixed Haar measure on $N$
(respectively on $A$). These are normalised so that under the
Iwasawa decomposition $ G=KAN$ we have $dg = e^{2\rho (\log a )}
dk da dn $ where $dk$ is the normalised Haar measure of $K$ and
$\log a$ is the unique element in ${\mathfrak a} $ such that
$\exp(\log a )=a$.

We follow the practice of using $C, C'$ etc. to denote constants
whose values are not necessarily the same at each occurrence.
\begin{Definition} For suitable  functions $f$on $X$, the Helgason
Fourier transform $\widetilde{f}$ of $f$   is defined by
$$\widetilde{f}(\lambda,
k)=\int_Xe^{(-i\lambda+\rho)(A(x,k))}f(x)dx; \lambda\in \mathfrak
a^*, k\in K.$$
\end{Definition}
Note that $\widetilde{f}$ descends to a function on $\mathfrak
a^*\times K/M$. By abuse of notation we will continue to denote
this function by $\widetilde{f}$. For $f\in L^1(X)$, there exists
a subset $B$ of $K$ of full Haar measure, such that
$\widetilde{f}(\lambda, k)$ exists for all $k\in B$ and
$\lambda\in \C$ with $|\Im \lambda|\le \rho$. Indeed for each
fixed $k\in B$, $\lambda\mapsto \widetilde{f}(\lambda, k)$ is
holomorphic in the strip $\{\lambda\in \C\,\, |\,\, |\Im
\lambda|<\rho\}$ and continuous on its boundary (see \cite{MRSS}).
\begin{Definition} For suitable functions $f$ on $X$, the Radon
transform ${\mathcal R}f$ of $f$ is defined by
$${\mathcal R}f(k, a)=e^{\rho(\log a)}\int_N f(kan)dn; k\in K, a\in A.$$
\end{Definition}
 ${\mathcal R}f$ descends
to a function on $K/M\times A$ and (as in the case of
$\widetilde{f}$) we continue to denote this function by ${\mathcal
R}f$. We will use the notation ${\mathcal R}f(k, t)$ for
${\mathcal R}f(k, a_t)$, $t\in \R$.

For a suitable function $f$ on $X$, the basic relation between
${\mathcal R}f$ and $\widetilde{f}$ is the following:
\begin{equation}\widetilde{f}(\lambda, k)={\mathcal
F}{\mathcal R}f(k, \cdot)(\lambda),
\label{Radon-to-Helgason}\end{equation} where ${\mathcal F}$
denotes the Euclidean Fourier transform on $A\equiv \R$.

Let $\what{K}_0$ be the set of equivalence classes of irreducible
unitary representations of $K$ which are class 1 with respect to
$M$, that is contains an $M$-fixed vector. Let $\delta\in
\what{K}_0$ and let $f\in L^1(X)$ be  $K$-finite of type $\delta$.
Then we have $d(\delta)\bar{\chi}_\delta*f=f$ where $d(\delta)$
(respectively $\chi_\delta$) denotes the degree (respectively
character) of $\delta$ and ($d(\delta)\bar{\chi}_\delta
*f)(x)=d(\delta)\int_Kf(kx)\bar{\chi}_\delta(k)dk$ for $x\in X$.
In particular if $\delta$ is the trivial representation then $f$
is a $K$-invariant function on $X$. For a function $f$ of type
$\delta$  we have $|f(x)|\le C\int_K|f(kx)|dk$ where
$C=d(\delta)\sup_{k\in K}|\chi_\delta(k)|\le d(\delta)^2$. Let
$g(x)=\int_K|f(kx)|dk$. Then $g\in L^1(X)$ and $g$ is
$K$-invariant, that is $g\in L^1(G)$ and $g$ is $K$-biinvariant.
We have $|f(x)|\le d(\delta)^2 g(x)$.
\begin{Definition}
For a $\delta$-type function $f$ in $L^1(X)$, the Abel transform
${\mathcal A}f$ of $f$ is defined by
$${\mathcal A}f(a)=e^{\rho(\log a)}\int_Nf(an)dn, \,a\in A.$$
\end{Definition}
It is well known that for a $K$-invariant function $g\in L^1(X)$,
${\mathcal A}g$ exists for almost every $a\in A$ and ${\mathcal
A}g\in L^1(A)$.  Now since
$$|{\mathcal A}f(a)|\le e^{\rho(\log a)}\int_N|f(an)|dn\le
d(\delta)^2e^{\rho(\log a)}\int_Ng(an)dn=d(\delta)^2{\mathcal
A}g(a)$$  for the $K$-invariant function $g$ constructed from $f$
as above, we conclude that  ${\mathcal A}f\in L^1(A)$. We will
write ${\mathcal A}f(t)$ for ${\mathcal A}f(a_t)$, $t\in \R$.

It is also well known that for $f\in L^1(X)$, ${\mathcal R}f\in
L^1(K\times A, dkda)$. We include the proof here for the sake of
completeness.   For $f\in L^1(X)$, we construct the function
$g(x)=\int_K|f(kx)|dk$. Then $g$ is a  $K$-biinvariant function in
$L^1(G)$ and hence as mentioned above (recall the identification
of $A$ and $\R$) ${\mathcal A}g(t)\in L^1(\R)$. But ${\mathcal
A}g(t)=e^{\rho t}\int_N g(a_tn)dn=e^{\rho
t}\int_N\int_K|f(ka_tn)|dkdn=\int_Ke^{\rho
t}\int_N|f(ka_tn)|dndk=\int_K{\mathcal R}|f|(k, t)dk$. Since
$\int_A{\mathcal A}g(t)dt<\infty$ we have $\int_A\int_K{\mathcal
R}|f|(k, t)dkdt<\infty$. This proves our assertion.

For $\lambda\in \mathfrak a^*_\C\equiv \C$, we denote by
$\phi_\lambda$ the elementary spherical function with parameter
$\lambda$. We have for all $x\in X$,
$\phi_\lambda(x)=\int_Ke^{(i\lambda+\rho)(A(x, k))}dk$ (see
\cite{He} p. 418). We will often regard $\phi_\lambda$ as a
$K$-biinvariant function on $G$.

For $x\in G$, we define $\sigma(x)=d(xK, K)$ where $d$ is the
canonical distance function for $X=G/K$ coming from the Riemannian
structure induced by the Cartan-Killing form restricted to
$\mathfrak p$. Here $\mathfrak g=\mathfrak k\oplus \mathfrak p$
(Cartan decomposition) and $\mathfrak p$ can be identified with
the tangent space at $eK$ of $G/K$.

The following estimates on the growth of $\phi_\lambda$ are well
known (\cite{He2}, \cite{GV} proposition 4.6.1 and theorems 4.6.4,
4.6.5): Let $\Im \lambda$ denote the imaginary part of $\lambda\in
\C$ and let $\Xi(x)=\phi_0(x)$. Then,
\begin{enumerate}
\item[(a)] $|\phi_\lambda(x)|\le 1$ for $\lambda\in \C, |\Im
\lambda|\le \rho$,

\item [(b)]$|\phi_\lambda(x)|\le e^{|\Im
\lambda|\sigma(x)}\Xi(x)$; for all $\lambda\in \C$

\item[(c)] $\Xi(a)\le C(1+\sigma(a))e^{-\rho(\log a)}$ for $a\in
\overline{A^+}=\exp \bar{\mathfrak a}_+$.
\end{enumerate}
We denote the  spherical Plancherel measure on $\mathfrak a^*$ by
$d\mu(\lambda)=\mu(\lambda)d\lambda$, where $d\lambda$ is the
Lebesgue measure. We have $\mu(\lambda)=|c(\lambda)|^{-2}$ where
$c(\cdot)$ is Harish-Chandra's $c$-function.

Recall that the  elements $\delta\in \what{K}_0$ are parametrized
by a pair of integers $(p_\delta, q_\delta)$ where $p_\delta\ge 0$
and $p_\delta\pm q_\delta\in 2\Z^+$ (see \cite{Kos, JW}). The
trivial representation in $\what{K}_0$ is parametrized by $(0,0)$
in this setup.

It is  known that for each $\delta\in \what{K}_0$, the $M$-fixed
vector is unique upto a scalar multiple (see \cite{Kos}). Let
$(\delta, V_\delta)\in \what{K}_0$. Suppose $\{v_i|i=1,\dots ,
d(\delta)\}$ is an orthonormal basis of $V_\delta$ of which $v_1$
is the $M$-fixed vector. Let $Y_{\delta,j}(k)=\langle
v_j,\delta(k)v_1\rangle, 1\le j\le d(\delta)$ and let $Y_0$ be the
$K$-fixed vector. Note that $Y_{\delta, j}$ is right $M$-invariant
that is it is a function on $K/M$. Recall that $L^2(K/M)$ is the
carrier space of the spherical principal series representations
$\pi_\lambda, \lambda\in \C$ in the compact picture and
$\{Y_{\delta,j}:1\le j\le d(\delta), \delta\in \what{K}_0\}$ is an
orthonormal basis for $L^2(K/M)$ adapted to the decomposition
$L^2(K/M)=\Si_{\delta\in \what{K}_0}V_\delta$. As the space $K/M$
can be identified with $S^{m_\gamma+m_{2\gamma}}$, this
decomposition can be viewed as the spherical harmonic
decomposition and therefore $Y_{\delta,j}$'s can be considered as
the {\em spherical harmonics}. The action of $\pi_\lambda$ is
given by:
$$(\pi_\lambda(x)g)(k)=e^{(i\lambda+\rho)A(x, k)}g(\kappa(x^{-1}k))\mbox{ for }
x\in G, k\in K \mbox{ and } g\in L^2(K/M).$$ Here $\kappa(x)$ is
the $K$-part of an element $x\in G$ in the Iwasawa decomposition
$G=KAN$. The representation  $\pi_\lambda$ is unitary for
$\lambda\in \R$. For $f\in L^1(X)$, $\delta\in \what{K}_0$ and
$1\le j\le d(\delta)$ we define,
$$f_{\delta, j}(x)=\int_Kf(kx)Y_{\delta, j}(k)dk.$$ It can be verified that
$f_{\delta, j}$ is a function of type $\delta$. The function $f$
can be  decomposed  as $f=\sum_{\delta\in
\what{K}_0}\sum_{j=1}^{d(\delta)}f_{\delta,j}$. In fact when $f\in
C^\infty(G)$  this is an absolutely convergent series in the
$C^\infty$-topology. When $f\in L^p(G), p\in [1,\infty)$,  the
equality is in the sense of distributions.

We have $|f_{\delta,j}(x)|\leq \|
Y_{\delta,j}\|_{\infty}\int_K|f(kx)|dk \le \int_K|f(kx)|dk$ as $\|
Y_{\delta,j}\|_{\infty}=\sup_{k\in K}|Y_{\delta,j}(k)|\le 1$.

For $\delta\in \what{K}_0$, $1\le j\le d(\delta)$,  $\lambda\in
\mathfrak a^*_\C$ and $x\in X$, we define
\begin{equation}\Phi_{\lambda,
\delta}^j(x)=\int_{K/M}e^{(i\lambda+\rho)(A(x,kM))}Y_{\delta,j}(kM)dk.
\label{Phi^j}\end{equation} We have, $\Phi_{\lambda, \delta}^j(x)=
\langle Y_{\delta,j}, \pi_{-\bar{\lambda}}(x)Y_0\rangle$, that is
$\Phi_{\lambda, \delta}^j$ is a matrix coefficient of the
spherical principal series in the compact picture. It is well
known (see \cite{He2}) that $\Phi^j_{\lambda, \delta}$'s are
eigenfunctions of the Laplace-Beltrami operator $\Delta$ with
eigenvalues $-(\lambda^2+\rho^2)$.  When $\delta=\delta_0$ is
trivial then $Y_{\delta_0,j}=Y_{\delta_0,1}=Y_0$ and
$\Phi_{\lambda, \delta_0}^1$ is obviously the elementary spherical
function $\phi_\lambda(x)$.
For $\lambda\in \mathfrak a^*_\C$, $x=ka_tK\in X$ and $1\le j\le
d(\delta)$, (see \cite{He2}, p. 344)
\begin{equation}\Phi^j_{\lambda,\delta}(x)
=Y_{\delta,j}(kM)\Phi^1_{\lambda,\delta}(a_t).
\label{helgason-to-jacobi}
\end{equation}
 $\Phi^1_{\lambda, \delta}$ is related with $\Phi^1_{-\lambda,
 \delta}$ by:
 \begin{equation}
\Phi^1_{\lambda,
\delta}=\frac{Q_\delta(\lambda)}{Q_\delta(-\lambda)}\Phi^1_{-\lambda,
 \delta},
 \label{intertwining}
 \end{equation} where $Q_\delta$'s are Kostant polynomials (see \cite{He2} p. 344,
 Theorem). Indeed $Q_\delta$ is the polynomial factor of $\Phi^1_{\lambda, \delta}$
 and hence of $\Phi^j_{\lambda,\delta}$ for $1\le j\le d(\delta)$ (see \cite{He2},
 p.344).
 The Kostant polynomial $Q_\delta$ is
 given by $$Q_\delta(\lambda)=(\frac 12(a+b+1+i\lambda))_{\frac{p_\delta+q_\delta}{2}}
 (\frac 12(a-b+1+i\lambda))_{\frac{p_\delta-q_\delta}{2}},$$
 where $(z)_m=\Gamma(z+m)/\Gamma(z)$,
 $a=\frac{m_\gamma+m_{2\gamma}-1}{2}$ and
 $b=\frac{m_{2\gamma}-1}{2}$. Thus $\deg Q_\delta=p_\delta$.

 Because of the relation (\ref{helgason-to-jacobi}) above we have
$\Phi^j_{\lambda,
\delta}=\frac{Q_\delta(\lambda)}{Q_\delta(-\lambda)}\Phi^j_{-\lambda,
 \delta}$.

We define the $j$-th $\delta$-spherical Fourier transform of $f$
by
$$\what{f}(\lambda)_{\delta,j}=\int_Xf(x)\Phi^j_{-\lambda,
\delta}(x)dx$$ for $\lambda \in {\mathfrak a^*}$. It is clear that
the $j$-th $\delta$-spherical Fourier transform of $f$ is the
$(\delta, j)$-th matrix coefficient of the operator Fourier
transform $\what{f}(\lambda)=\int_Gf(x)\pi_{-\lambda}(x)dx$. It is
not difficult to verify that
$$\what{f}(\lambda)_{\delta,j}=\int_Xf_{\delta, j}(x)\Phi^j_{-\lambda,
\delta}(x)dx=\what{f}_{\delta, j}(\lambda) \mbox{ say. Furthermore
} \int_Xf_{\delta', j'}(x)\Phi^j_{-\lambda, \delta}(x)dx=0,$$ when
$\delta\neq \delta'$ or $j\neq j'$. Henceforth we will not
distinguish between $\what{f}(\lambda)_{\delta,j}$ and
$\what{f}_{\delta, j}(\lambda)$. Let $f\in L^1(X)\cap L^2(X)$.
Then for almost every $\lambda\in \mathfrak a^*$,
\begin{equation}
\|\what{f}(\lambda)\|_2^2=\sum_{\delta\in \what{K}_0}\sum_{1\le
j\le d(\delta)}|\what{f}_{\delta, j}(\lambda)|^2,
\label{decompose-HS-norm}
\end{equation}
where $\|\,\cdot\,\|_2$ is the Hilbert-Schmidt norm.

Also by (\ref{intertwining})
$$\what{f}_{\delta,j}(-\lambda)=\frac{Q_\delta(\lambda)}{Q_\delta(-\lambda)}
\what{f}_{\delta,j}(\lambda).$$ Note that (see \cite{He2} p. 348)
for $\lambda\in\mathfrak a^*$,
$\overline{Q_\delta(\lambda)}=Q_\delta(-\lambda)$. Consequently
$|\what{f}_{\delta, j}(\lambda)|=|\what{f}_{\delta,
j}(-\lambda)|$.

The following is also easy to see:
\begin{equation}\ba{ll}&\int_K \widetilde{f}(\lambda,
k)Y_{\delta,j}(k)dk\\ \\
=&\int_X\int_Kf(x)e^{(-i\lambda+\rho)A(x,k)}Y_{\delta,j}(k)dkdx\mbox{ (by Fubini's theorem)}\\
\\
=&\int_Xf(x)\Phi^j_{-\lambda,\delta}(x)dx\\ \\
=&\what{f}_{\delta,j}(\lambda).\ea
\label{Helgason-to-delta-spherical}
\end{equation}

Starting from the relation (\ref{Radon-to-Helgason}) and  using
(\ref{Helgason-to-delta-spherical}) we have,
$$\int_K{\mathcal F}({\mathcal
R}(f)(k,\cdot))(\lambda)Y_{\delta,j}(k)dk=\what{f}_{\delta,j}(\lambda).$$
Now the left hand side is (recall that $\mathfrak a^*\equiv \R$):
\begin{equation}\ba{ll}
&\int_K\int_R{\mathcal R}(f)(k,t)e^{-i\lambda t}dt\,
Y_{\delta,j}(k)dk\\ \\
=&\int_K\int_R e^{\rho t}\int_Nf(ka_tn)dne^{-i\lambda t}dt\,
Y_{\delta, j}(k)dk\\ \\
=&\int_R e^{\rho t}\int_Nf_{\delta,j}(a_tn)dn\,e^{-i\lambda t}dt\mbox{ (by Fubini's theorem)}\\
\\
=&\int_\R{\mathcal A}(f_{\delta,j})(t)e^{-i\lambda t}dt\\ \\
=&{\mathcal F}({\mathcal A}(f_{\delta,j}))(\lambda). \ea$$
Therefore $${\mathcal F}({\mathcal
A}(f_{\delta,j})(\lambda)=\what{f}_{\delta,j}(\lambda).
\label{abel-to-delta-j}
\end{equation}

Note that from above it is also clear that:
\begin{equation}\int_K {\mathcal
R}(f)(k, t)Y_{\delta, j}(k)dk={\mathcal A}(f_{\delta, j})(t)
\label{radon-to-abel}
\end{equation}
and hence
\begin{equation}
|{\mathcal A}(f_{\delta, j})(t)|=\left |\int_K{\mathcal R}(f)(k,
t)Y_{\delta, j}(k)dk\right |\le \int_K |{\mathcal R}(f)(k,
t)|dk\le \int_K {\mathcal R}|f|(k, t)dk \label{compare-radon-abel}
\end{equation} since $\|Y_{\delta, j}\|_\infty\le 1$.

We will conclude this section with  a description of the {\em
heat-kernel} of the symmetric space $X$. The heat kernel on $X$ is
an appropriate analogue of the Gauss kernel $p_t$
 on $\R^n$ where $p_t(x)=(4\pi t)^{-\frac{n}{2}}e^{-\frac{\|x\|^2}{4t}}, t>0$.

Let $\Delta$ be the Laplace-Beltrami operator of $X$. Then  (see
\cite{St}, Chapter V), $T_t=e^{t\Delta},t>0$ defines a semigroup
(heat-diffusion semigroup) of operators such that for any $\phi\in
C^\iy_c(X)$, $T_t\phi$ is a solution of $\Delta u=\f{\partial u
}{\partial t}$ and $T_t\phi\lgra \phi$ a.e. as $t\lgra 0$. For
every $t>0$, $T_t$ is an integral operator with kernel $h_t$, that
is for any $\phi\in C^\iy_c(X)$, $T_t\phi=\phi*h_t$. The $h_t,
t>0$ are $K$-biinvariant functions on $G$, $h(x,t)=h_t(x)$ as a
function of the variables $t\in \R^+$ and $x\in G$ is in
$C^\iy(G\times \R^+)$ and has the following  properties:
\begin{enumerate}
\item[i.] $\{h_t: t>0\}$ form a semigroup under convolution $*$.
That is $h_t*h_s= h_{t+s}$ for $t,s>0$. \item[ii.] $h_t$ is a
fundamental solution of $\Delta u=\f{\partial u}{\partial t}$.
\item[iii.] $h_t\in L^1(G)\cap L^\iy(G)$ for every  $t>0$.
\item[iv.] $\int_{X} h_t(x) dx =1$ for every $t>0$.
\end{enumerate}

Thus we see that the heat kernel $h_t$ on $X$ retains all the nice
properties of the Gauss  kernel.  It is well known that $h_t$ is
given by (see e.g. \cite{A4}):
\begin{equation}\label{heatkernel}
h_t(x) = \frac{1}{|W|}\int_{\mathfrak
a^*}e^{-t(\lambda^2+\rho^2)}\phi_\lambda(x)\mu(\lambda)d\lambda.
\end{equation}
 That is, the spherical Fourier transform of $h_t$,
$\what{h_t}(\lambda)= e^{-t(\lambda^2+\rho^2)}$. It has been
proved in \cite{A4} (Theorem 3.1 (i)) that for any $t>0$, there
exists $C>0$ depending only on $X$ such that
\begin{equation}\label{ankerestimate}h_t(\exp H)
 \leq Ct^{-\f{1}{2}} e^{-\rho^2t-\langle\rho,
H\rangle-\f{|H|^2}{4t}}(1+|H|^2)^{\f{d_X-1}{2}}\end{equation} for
 $H\in
 \overline{\mathfrak a^+}$, where $d_X=m_\gamma+m_{2\gamma}+1=\dim X$.

\section{statement and proof of the theorem}\newsection
\begin{Theorem}
Let $f\in  L^2(X)$ satisfy
\begin{equation}\int_X\int_{\mathfrak a^*}
\frac{|f(x)|\|\what{f}(\lambda)\|_2e^{\sigma(x)|\lambda|}\Xi(x)}{(1+\sigma(x)+|\lambda|)^d}
dxd\mu(\lambda)<\infty \label{beurling-heat-condition}
\end{equation} for some nonnegative integer $d$.
Then $f$ is a  $K$-finite function of the form $f=\sum_{\delta\in
F}h_\delta$ where $F=\{\delta\in \what{K}_0\,|\, p_\delta<\frac
{d-d_X} 2\}$ is a finite set of $K$-types, $h_\delta$ is a
function of type $\delta$ having Fourier coefficients
$\what{h}_{\delta,
j}(\lambda)=P'_{\delta,j}(\lambda^2)Q_\delta(\lambda)e^{-\alpha\lambda^2}$
for $1\le j\le d(\delta)$. Here $\alpha$ is a positive constant
and $P'_{\delta,j}$ a polynomial which  depends on $\delta$ and
$j$. Consequently $f$ is a derivative of the heat kernel
$h_\alpha$.

In particular if  $d\le d_X$, then $f= 0$ almost everywhere.
\label{thm-symm-heat}
\end{Theorem}

\begin{proof}
We have divided the proof in several steps for the convenience of
the readers. We will use Fubini's theorem freely throughout the
proof without explicitly mentioning it.

\noindent {\bf Step 1:} In this step we will show that $f\in
L^1(X)$.

Note that, since $f\in L^2(X)$, $f$ is a locally integrable
function on $X$. Now there are two possible situations.
\begin{enumerate}
\item[(a)] $\what{f}$ is  supported on a set of infinite measure.

\item[(b)] $\what{f}$ is  supported on  a set of finite measure.
\end{enumerate}
In case (a) clearly $f\in L^1(X)$, since from
(\ref{beurling-heat-condition}) it follows that there exists
$\lambda\in \mathfrak a^*$, $|\lambda|>2\rho$ such that $\int_X
|f(x)|e^{(|\lambda|-\rho)\sigma(x)}dx<\infty$.

In case (b), as $\what{f}\not\equiv 0$, there exists
$\lambda_0\neq 0$ such that $\what{f}(\lambda_0)\neq 0$. Suppose
$|\lambda_0|=r>0$. Then from (\ref{beurling-heat-condition})  we
have $\int_X\frac{|f(x)|e^{r\sigma(x)}\Xi(x)}{(1+\sigma(x)+r)^d}
dx<\infty$. As $f\in L^1_{loc}(X)$, for $0<r'<r$, $\int_X
|f(x)|e^{r'\sigma(x)}\Xi(x)dx<\infty$.  For  $\delta\in
\what{K}_0$ and $1\le j\le d(\delta)$, we consider the integral
$\what{f}_{\delta,j}=\int_Xf(x) \Phi_{-\lambda, \delta}^j(x)dx$.
Then
$$\ba{lll}|\int_Xf(x) \Phi_{-\lambda, \delta}^j(x)dx|&\le& \int_X|f(x)||\Phi_{-\lambda,\delta}^j(x)|dx\\ \\
&\le&\int_X|f(x)|e^{|\Im \lambda|\sigma(x)}\Xi(x)\\ \\
&\le &\int_X|f(x)|e^{r'\sigma(x)}\Xi(x)e^{(|\Im
\lambda|-r')\sigma(x)}dx.\ea$$

This shows that $\what{f}_{\delta,j}$ is analytic in the open
strip $|\Im \lambda|<r'$ in $\mathfrak a^*_\C$, which
 contradicts the assumption that $\what{f}$ and hence $\what{f}_{\delta,j}$ is
 supported on a set of finite measure. This completes step 1.

 \noindent{\bf Step 2:} In this step we will show that (\ref{beurling-heat-condition}) is equivalent to the condition:
\begin{equation}I=\int_G\int_{\mathfrak a^*}
\frac{|f(x)||\what{f}_{\delta,j}(\lambda)|e^{\sigma(x)|\lambda|}
e^{-\rho H(x)}}{(1+\sigma(x)+|\lambda|)^d}dxd\mu(\lambda)<\infty
\label{starting-point}
\end{equation} for  every  fixed $\delta\in \what{K}_0$
and $1\le j\le d(\delta)$.

 We know that
$\|\what{f}(\lambda)\|_2^2=\sum_{\delta\in
\what{K}_0}\sum_{j=1}^{d(\delta)}|\what{f}_{\delta,
j}(\lambda)|^2$.  Therefore we obtain the following from
(\ref{beurling-heat-condition})
$$I_1=\int_G\int_{\mathfrak a^*}
\frac{|f(x)||\what{f}_{\delta,j}(\lambda)|e^{\sigma(x)|\lambda|}
\Xi(x)}{(1+\sigma(x)+|\lambda|)^d}dxd\mu(\lambda)<\infty.$$ We
have changed the integration over $X$ to integration over $G$ as
all the terms of the integrand are right $K$-invariant.
 As
 $f$ and $\sigma$  are both right
$K$-invariant, on replacing $x$ by $xk^{-1}$ in the integral $I$
we  get
$$I=I(k)=\int_G\int_{\mathfrak a^*}
\frac{|f(x)||\what{f}_{\delta,j}(\lambda)|e^{\sigma(x)|\lambda|}
e^{-\rho H(xk^{-1})}}{(1+\sigma(x)+|\lambda|)^d}dxd\mu(\lambda).$$
Therefore $$I=\int_K I(k)dk=\int_K \int_G\int_{\mathfrak a^*}
\frac{|f(x)||\what{f}_{\delta,j}(\lambda)|e^{\sigma(x)|\lambda|}
e^{-\rho
H(xk^{-1})}}{(1+\sigma(x)+|\lambda|)^d}dxd\mu(\lambda)dk.$$ Since
 $\int_Ke^{-\rho
H(xk^{-1})}dk=\Xi(x^{-1})$ and  $\Xi(x)=\Xi(x^{-1})$, we conclude
that $I=I_1<\infty$. \vspace{.25in}

\noindent{\bf Step 3:} We will now show that:
\begin{equation}
\int_K\int_\R\int_{\R}\frac{{\mathcal
R}(|f|)(k,t)|\what{f}_{\delta,j}(\lambda)|e^{|\lambda||t|}}{(1+|t|+|\lambda|)^d}dkdtd\mu(\lambda)
<\infty. \label{target-heat0}
\end{equation}
Since  the  integrand is even in $\lambda$ as pointed out before,
this is equivalent to showing
\begin{equation}
\int_K\int_{\R}\int_{\R^+}\frac{{\mathcal
R}(|f|)(k,t)|\what{f}_{\delta,j}(\lambda)|e^{|\lambda||t|}}{(1+|t|+|\lambda|)^d}dkdtd\mu(\lambda)\,<\infty
\label{target-heat1}
\end{equation}  We will break the above
integral  into the following 3 parts and show that each part is
finite. That is we will show:
\begin{enumerate}
\item[(i)]
$$\int_K\int_{\R}\int_{L}^{\infty}\frac{{\mathcal
R}(|f|)(k,t)|\what{f}_{\delta,j}(\lambda)|e^{|\lambda||t|}}
{(1+|t|+|\lambda|)^d}dtd\mu(\lambda)dk<\infty$$ for $L>0$ such
that $L^2+L>d$.

\item[(ii)]
$$\int_K\int_{|t|>M}\int_{0}^{L}
\frac{{\mathcal
R}(|f|)(k,t)|\what{f}_{\delta,j}(\lambda)|e^{|\lambda||t|}}
{(1+|t|+|\lambda|)^d}dtd\mu(\lambda)dk<\infty$$ for
$M=2(L+1+\rho)$ and $L$ as in (i).

\item[(iii)]
$$\int_K\int_{|t|\leq M}\int_0^{L}\frac{{\mathcal
R}(|f|)(k,t)|\what{f}_{\delta,j}(\lambda)|
e^{|\lambda||t|}}{(1+|t|+|\lambda|)^d}dtd\mu(\lambda)dk<\infty$$
for $M,L$ used in (i) and (ii).
\end{enumerate}  If we show that  (i),
(ii) and  (iii) are finite, they together will obviously imply
(\ref{target-heat1}) and hence (\ref{target-heat0}).
\vspace{.08in}

\noindent{\em Proof of} (iii): As the domain of integration
$[-M,M]\times [0, L]$ is compact and as $\frac {|\what{f}_{\delta,
j}(\lambda)|e^{|\lambda||t|}}{(1+|\lambda|+|t|)^d}$ is continuous
in this domain, the integral is bounded by
$C\int_K\int_{-M}^M{\mathcal R}(|f|)(k,t)dt dk$. Recall that $f\in
L^1(G)$. Therefore,
$$\ba{lll}\int_K\int_\R {\mathcal R}(|f|)(k,t)e^{\rho
t}dtdk&=&\int_K\int_\R e^{\rho t}\int_N |f(ka_tn)|dne^{\rho
t}dtdk\\&=&\int_K\int_\R\int_N |f(ka_tn)|e^{2\rho t}dkdtdn
\\&=&\int_G|f(g)|dg<\infty.\ea$$ Hence
$\int_K\int_{-M}^M{\mathcal R}(|f|)(k,t)dtdk<\infty$, since this
is $\le 2e^{\rho M}\int_K\int_\R{\mathcal R}(|f|)(k,t)e^{\rho t}dt
<\infty$. \vspace{.08in}

\noindent{\em Proof of} (i): It is given that $L+L^2>d$. We will
show that for any $\lambda$ such that $|\lambda|\geq L$,
\begin{equation}\frac{e^{|\lambda|\sigma(a_tn)}}{(1+|\lambda|+\sigma(a_tn))^d}\geq
\frac{e^{|\lambda|\sigma(a_t)}}{(1+|\lambda|+\sigma(a_t))^d}.
\label{useful-inequality}
\end{equation}
 Let
$F(x)=\frac{e^{\alpha x}}{(1+\alpha+x)^d}$ for $\alpha>0$ and
$\alpha+\alpha^2>d$. Then $F'(x)>0$ for any $x\geq 0$. If $x\geq
y\geq 0$, then
\begin{equation}\frac{e^{\alpha
x}}{(1+\alpha+x)^d}\geq \frac{e^{\alpha y}}{(1+\alpha+y)^d}.
\label{elementary}
\end{equation}

Note that $\sigma(an)\geq \sigma(a)$ for all $a\in A$ and $n\in
N$. Now take $x=\sigma(a_tn)$ and $y=\sigma(a_t)$. Then $x\geq
y\geq 0$. We take $\alpha=|\lambda|\geq L$ to get the required
result.

We start now from (\ref{starting-point}) and use the Iwasawa
decomposition $G=KAN$ and the inequality (\ref{useful-inequality})
to obtain:
\begin{equation}\int_K\int_\R\int_{L}^\infty\frac{{\mathcal
R}(|f|)(k,t)|\what{f}_{\delta,j}(\lambda)|e^{|\lambda|t}}
{(1+|t|+|\lambda|)^d}dtd\mu(\lambda)dk<\infty.
\label{required-later}
\end{equation}
This proves (i). \vspace{.08in}

\noindent{\em Proof of} (ii): Let
$$I_2=\int_K\int_{|t|>M}\int_0^{L}\frac{{\mathcal
R}(|f|)(k,t)|\what{f}_{\delta,
j}(\lambda)|e^{|\lambda||t|}}{(1+|t|+|\lambda|)^d}dkdtd\mu(\lambda).$$
Since $|\what{f}_{\delta, j}(\lambda)|$ is bounded and
$\mu(\lambda)$ is continuous, $I_2\leq C\int_K\int_{|t|>M}
\frac{{\mathcal
R}(|f|)(k,t)|e^{L|t|}}{(1+|t|)^d}dkdt=C\int_K\int_{|t|>M} \int_N
\frac{|f(ka_tn)|e^{L|t|}e^{\rho t}}{(1+|t|)^d}dkdtdn=CI_3$, say.
We will show that $I_3$ is finite for  $M=2(L+1+\rho)$.

We start  with the assumption that $\what{f}_{\delta, j}\not\equiv
0$ (otherwise  $f_{\delta, j}= 0$ almost everywhere). This implies
that $\what{f}_{\delta, j}(\lambda)\neq 0$ for almost every
$\lambda\in \R$ as $\what{f}_{\delta, j}$ is real analytic on
$\R$, $f$ being an $L^1$-function.

Therefore using from (\ref{starting-point}) we can get a
$\lambda_0\in \R$ with $|\lambda_0|>2(L+\rho)$ such that:
$$\int_G\frac{|f(x)|e^{\sigma(x)|\lambda_0|}e^{-\rho H(x)}}{(1+\sigma(x)+|\lambda_0|)^d}<\infty.$$

We will use the Iwasawa decomposition $G=KAN$. Now since $\sigma$
and $H$ are both left $K$-invariant and $H$ is right $N$-invariant
we  obtain:
$$\int_K\int_{\R\times N}\frac{|f(ka_tn)|e^{|\lambda_0|\sigma(a_tn)}e^{-\rho t}}
{(1+\sigma(a_tn)+|\lambda_0|)^d}e^{2\rho t}dkdtdn<\infty.$$ Notice
that $|\lambda_0|+|\lambda_0|^2>d$. Therefore by applying the
argument of case (i) (see (\ref{elementary})) to $|\lambda_0|$ we
get:
$$\frac{e^{|\lambda_0|\sigma(a_tn)}}{(1+|\lambda_0|+\sigma(a_tn))^d}\geq
\frac{e^{|\lambda_0|\sigma(a_t)}}{(1+|\lambda_0|+\sigma(a_t))^d}.$$

Therefore in particular:
$$\int_K\int_{|t|>M}\int_{N}\frac{|f(ka_tn)|
e^{|\lambda_0||t|}e^{\rho
t}}{(1+|t|+|\lambda_0|)^d}dkdtdn<\infty.$$ Note that $M+M^2>d$ as
$M=2(L+1+\rho)$ and $L+L^2>d$. Applying the argument of case (i)
again (see (\ref{elementary})) this time with  $\alpha=|t|>M$ and
$x=|\lambda_0|, y=2(L+\rho)$ we get,
$$\frac{e^{|\lambda_0||t|}}{(1+|\lambda_0|+|t|)^d}\geq
\frac{e^{2(L+\rho)|t|}}{(1+2(L+\rho)+|t|)^d}.$$ Therefore,
$$\int_K\int_{|t|>M}\int_{N}
\frac{|f(ka_tn)|e^{2(L+\rho)|t|}e^{\rho
t}}{(1+|t|+2(L+\rho))^d}dkdtdn<\infty.$$
We also see  that for $|t|>M=2(L+1+\rho)$
$$\frac{e^{(L+\rho)|t|}}{(1+|t|+2(L+\rho))^d}>\frac 1{(1+|t|)^d}$$
(This  is equivalent to showing: $e^{(L+\rho)|t|}>(1+\frac
{2(L+\rho)}{1+|t|})^d$. For $|t|> M=2(L+1+\rho)$,
$e^{(L+\rho)|t|}>e^{2(L+\rho)(L+1+\rho)}>e^d>2^d$ (as $L+L^2>d$),
while $(1+\frac
{2(L+\rho)}{1+|t|})^d<(1+\frac{2(L+\rho)}{1+2(L+1+\rho)})^d
<2^d$.)

So we obtain:
$$\int_K\int_{|t|>M}\int_{N}\frac{|f(ka_tn)|e^{(L+\rho)|t|}e^{\rho t}}{(1+|t|)^d}dk
dtdn<\infty$$ and hence, $I_3<\infty$. This completes the proof of
(ii).

Thus from (i), (ii) and (iii) we obtain (\ref{target-heat0}).
\vspace{.25in}

\noindent{\bf Step 4:} From (\ref{target-heat0}) and
(\ref{compare-radon-abel}) we have,

\begin{equation}
\int_\R\int_{\R}\frac{|{\mathcal A}(f_{\delta',
j'})(t)||\what{f}_{\delta,j}(\lambda)|e^{|\lambda||t|}}{(1+|t|+|\lambda|)^d}dtd\mu(\lambda)
<\infty \label{target-intermediate}
\end{equation}
for $\delta, \delta'\in \what{K}_0$, $1\le j\le d(\delta)$, $1\le
j'\le d(\delta')$.

In particular we can take $\delta=\delta'$ and $j=j'$ to obtain:
\begin{equation}
\int_\R\int_{\R}\frac{|{\mathcal A}(f_{\delta,
j})(t)||\what{f}_{\delta,j}(\lambda)|e^{|\lambda||t|}}{(1+|t|+|\lambda|)^d}dtd\mu(\lambda)
<\infty. \label{target-last}
\end{equation}
\vspace{.25in}

\noindent{\bf Step 5:} Now  we will show that in
(\ref{target-last}) $d\mu(\lambda)$ can be replaced by $d\lambda$.
We have the following asymptotic estimate of the spherical
Plancherel density (see \cite{Ank1})
\begin{equation}\label{estimate-plancherel-measure}
\mu(\lambda)=|c(\lambda)|^{-2}\asymp
\langle\lambda,\gamma\rangle^2
(1+|\langle\lambda,\gamma\rangle|)^{m_\gamma +m_{2\gamma}-2},
\end{equation} where $m_\gamma, m_{2\gamma}$ are as defined in
section 2.  Here $f\asymp g$ means $c_1g(\lambda)\le f(\lambda)\le
c_2g(\lambda)$ for two positive constants $c_1, c_2$ and
$\lambda\in\mathfrak a^*$, $|\lambda|$ large.

For some suitable large $R>0$ let
$$M_1(\lambda) =
\int_{|t|>R}\frac{|{\mathcal A}(f_{\delta, j})(t)|
 e^{|\lambda||t|}}{(1+|t|+|\lambda|)^d}dt,\,\,
M_2(\lambda) = \int_{|t|\le R}\frac{|{\mathcal A}(f_{\delta,
j})(t)|
 e^{|\lambda||t|}}{(1+|t|+|\lambda|)^d}dt $$ and let
 $M(\lambda)=M_1(\lambda)+M_2(\lambda)$.

Then we are given that $\int_{\mathfrak a^*}
M(\lambda)|\what{f}_{\delta,j}(\lambda)| |c(\lambda)|^{-2}
d\lambda $ is finite  which implies that the integrand is finite
for almost every $\lambda$. Now both $M_1(\lambda)$ and
$M_2(\lambda)$ are clearly radial. Note also that $M_1$ is an
increasing function of $|\lambda|$ (see (\ref{elementary})) and
$M_2(\lambda)\le e^{|\lambda|R}\times \|{\mathcal A}(f_{\delta,
j})\|_1$ for any $\lambda$. Therefore $M_2(\lambda)$ is bounded on
compact sets. The Plancherel density $|c(\lambda)|^{-2}$ is real
analytic. Furthermore since $f\in L^1(X)$ the function $\lambda
\rightarrow \what{f}_{\delta,j}(\lambda)$ is real analytic. Hence
the set of zeros of these functions is at most countable, in
particular they have measure zero. Therefore $M_1(\lambda)$ is
finite everywhere and locally integrable since it is an increasing
function of $|\lambda|$. Thus both $M_1$ and $M_2$ are locally
integrable and hence $M$ is locally integrable. We want to show
that $\int_{\mathfrak a^*}
M(\lambda)|\what{f}_{\delta,j}(\lambda)| d\lambda $ is finite. The
integrand is locally integrable on ${\mathfrak a^*}$, hence we
need only examine its behaviour for large $|\lambda|$. Now the
above-mentioned formula for the Plancherel density shows that $
|c(\lambda)|^{-2}$ tends to $\infty$ as $|\lambda |$ tends to
$\infty$. In particular there exists $ A > 0$ such that
$|c(\lambda)|^{-2}> 1$ whenever $|\lambda| > A $. Now $\infty
>\int_{|\lambda| > A} M(\lambda)|\what{f}_{\delta,j}(\lambda)|
|c(\lambda)|^{-2} d\lambda \geq \int_{|\lambda| > A} M(\lambda)
|\what{f}_{\delta,j}(\lambda)| d\lambda$. This immediately implies
our assertion, that is we get

\begin{equation}
\int_\R\int_{\R}\frac{|{\mathcal A}(f_{\delta,
j})(t)||\what{f}_{\delta,j}(\lambda)|e^{|\lambda||t|}}{(1+|t|+|\lambda|)^d}dtd\lambda
<\infty. \label{target-last-last}
\end{equation}
\vspace{.08in}

\noindent{\bf Step 6:} In this step we will deduce that
$\what{f}_{\delta,j}(\lambda)=P(\lambda)e^{-\alpha\lambda^2}$,
where $P$ is a polynomial which depends on $\delta, j$
 and $\alpha$ is a positive constant, which is independent of $\delta,j$.

Note that ${\mathcal A}(f_{\delta,j})\in L^1(\R)\cap L^2(\R)$. In
view of (\ref{abel-to-delta-j}) we can apply  theorem
\ref{thm-BDJ} to obtain
$\what{f}_{\delta,j}(\lambda)=P(\lambda)e^{-\alpha\lambda^2}$. A
priori the polynomial $P$ as well as the constant $\alpha$ depend
on $\delta, j$. We will see that the constant $\alpha$ is actually
independent of $\delta,j$.

Suppose if possible for $\delta_1, \delta_2\in \what{K}_0$ and
$1\le j_1\le d(\delta_1), 1\le j_2\le d(\delta_2)$,
$$\ba{ll}(1)&
\what{f}_{\delta_1,j_1}(\lambda)=P_1(\lambda)e^{-\alpha_1\lambda^2}\\&\\
(2)&\what{f}_{\delta_2,j_2}(\lambda)=P_2(\lambda)e^{-\alpha_2\lambda^2},\ea$$
where $P_1, P_2$ are two polynomials, $\alpha_1, \alpha_2$ are
positive constants and $\alpha_1\neq \alpha_2$. Without loss of
generality we can assume that $\alpha_1<\alpha_2$. From (2) above
we have,

$$\ba{ll}(3)& {\mathcal A}(f_{\delta_2,j_2})(t)=P_2(t)e^{-\frac
1{4\alpha_2}t^2}.\ea$$ Substituting (1) and (3) in
(\ref{target-intermediate}) we see that the integrand in
(\ref{target-intermediate}) is
$$\frac{|P_1(\lambda)||P_2(t)|e^{-(\sqrt{\alpha_1}|\lambda|-\frac 1{2\sqrt{\alpha_2}}|t|)^2}
e^{A|\lambda||t|}}{(1+|t|+|\lambda|)^d}$$ where
$A=1-\sqrt{\frac{\alpha_1}{\alpha_2}}>0$ as
$\frac{\alpha_1}{\alpha_2}<1$. Therefore the integrand in
(\ref{target-intermediate}) grows very rapidly in the
neighbourhood of the hyperplane (pair of straight lines)
$\sqrt{\alpha_1}|\lambda|=\frac 1{2\sqrt{\alpha_2}}|t|$ and the
integral diverges. This establishes that the positive constant
$\alpha$ is independent of $\delta$ and $j$. \vspace{.25in}

\noindent{\bf Step 7:} This is our final step wherein we conclude
the proof of the theorem. From the previous step we know that
$\what{f}_{\delta,j}(\lambda)=P(\lambda)e^{-\alpha\lambda^2}$.
This shows that  $f_{\delta,j}$ is a derivative of the heat kernel
$h_\alpha$. Notice also that
$P(\lambda)=P'(\lambda^2)Q_\delta(-\lambda)$ where $P'(\lambda^2)$
is a polynomial
 in $\lambda^2$, because $Q_\delta(-\lambda)$ is a factor of $\what{f}_{\delta,j}(\lambda)$ (see section 2).
 Recalling that $\deg Q_\delta=p_\delta$
 we see that $\deg P(\lambda)\ge p_\delta$.

On the other hand noting that ${\mathcal
A}(f_{\delta,j})(t)=P(t)e^{-\frac 1{4\alpha}t^2}$, substituting
$\what{f}_{\delta,j}$ and  ${\mathcal A}(f_{\delta,j})$ back in
(\ref{target-last}) and using (\ref{estimate-plancherel-measure})
it is easy to verify that  $\deg P<\frac {d'-1}2$ where
$d'=d-(m_\gamma+m_{2\gamma})$ as otherwise the integral in
(\ref{target-last}) diverges.

 Therefore if
 $p_\delta\ge \frac {d'-1}2$, then $f_{\delta, j}= 0$ almost everywhere.
 As $p_\delta\ge |q_\delta|$, we conclude that only for finitely
 many $\delta\in \what{K}_0$ which are parametrized by $(p_\delta, q_\delta)$ with
 $p_\delta<\frac {d'-1}2$,  $f_{\delta, j}$ will satisfy
 (\ref{target-last}). We thus conclude that $f$ is a $K$-finite function whose each $K$-isotypical component
 is a derivative of the heat kernel $h_\alpha$. In other words,
 $f$ itself is a derivative of the heat kernel $h_\alpha$.

 In particular if $d'\le 1$ that is if $d\le 1+m_\gamma+m_{2\gamma}=d_X$,
 then there is no $p_\delta$ satisfying $p_\delta<\frac {d'-1} 2$
 and hence in that case $f=0$ almost everywhere.
\end{proof}
\section{Sharpness of the estimate}\newsection
In order to complete the picture we investigate the optimality of
the condition used in theorem 3.1. More precisely, suppose a
function $f\in L^1(X)\cap L^2(X)$ satisfies
\begin{equation}\int_X\int_{\mathfrak a^*}
\frac{|f(x)\|\what{f}(\lambda)\|_2e^{c\sigma(x)|\lambda|}\Xi(x)^{1-\varepsilon}}{(1+\sigma(x)+|\lambda|)^d}
dxd\mu(\lambda)<\infty \label{beurling-heat-condition-sharp1}
\end{equation} for some nonnegative integer $d$ and $c, \varepsilon \in
\R$. Then:

\begin{enumerate}
\item[(i)] We will see that if $\{c>1 \mbox{ and } \varepsilon \ge
0\}$ or if $\{c\ge 1\mbox{ and }\varepsilon
>0\}$ in (\ref{beurling-heat-condition-sharp1}) then $f= 0$ almost everywhere.

\item[(ii)] We will find a symmetric space $X$ on which  there can
be infinitely many linearly independent functions in $L^1(X)\cap
L^2(X)$ satisfying the estimate
(\ref{beurling-heat-condition-sharp1}) with $\{c<1 \mbox{ and
}\varepsilon \le 0\}$ and with $\{c\le 1 \mbox{  and } \varepsilon
<0\}$. These functions are not of the form characterized in
theorem 3.1.
\end{enumerate}
In case (i) as $c>1$ and $\Xi^{-\varepsilon}\ge 1$, $f$ satisfies
the condition (\ref{beurling-heat-condition}) in theorem 3.1 and
hence $\what{f}_{\delta, j}(\lambda)=P_{\delta, j}(\lambda)
e^{-\alpha\lambda^2}$. Therefore ${\mathcal
A}(f_{\delta,j})(t)=P_{\delta,j}(t)e^{-\beta t^2}$ where
$\alpha\beta=\frac 14$, since ${\mathcal A}(f_{\delta,j})$ is the
Euclidean Fourier inverse of $\what{f}_{\delta, j}$.

On the other hand starting from the condition
(\ref{beurling-heat-condition-sharp1}) and following the steps of
the proof of theorem 3.1 we obtain finally:

\begin{equation}
\int_\R\int_{\R}\frac{|{\mathcal A}(f_{\delta,
j})(t)||\what{f}_{\delta,j}(\lambda)|e^{c|\lambda||t|}e^{\varepsilon
\rho t}}{(1+|t|+|\lambda|)^d}dtd\lambda <\infty.
\label{beurling-heat-condition-sharp1A}
\end{equation}
Substituting ${\mathcal A}(f_{\delta, j})$ and $\what{f}_{\delta,
j}$ as obtained above in this inequality we see that it demands
$$\int_\R\int_{\R}\frac{e^{-(\sqrt \alpha |\lambda|-\sqrt \beta|t|)^2}e^{(c-1)|\lambda||t|}e^{\varepsilon
\rho t}}{(1+|t|+|\lambda|)^d}dtd\lambda <\infty.$$ But around the
hyperplane $\sqrt \alpha |\lambda|=\sqrt \beta|t|$ the integrand
grows rapidly as $|t|\lgra \infty$ since $c-1>0$ or $\varepsilon
>0$ and hence the integral becomes infinite which contradicts (\ref{beurling-heat-condition-sharp1A}).

 Next we consider the case (ii) that is, we will find a symmetric space $X$ and functions  $f$ on $X$ which  satisfy
\begin{equation}\int_X\int_{\mathfrak a^*}
\frac{|f(x)\|\what{f}(\lambda)\|_2e^{c\sigma(x)|\lambda|}\Xi(x)^{1+\varepsilon'}}{(1+\sigma(x)+|\lambda|)^d}
dxd\mu(\lambda)<\infty \label{beurling-heat-condition-sharp2}
\end{equation} for some nonnegative integer $d$, and either $c<1, \varepsilon' \ge 0$ or $c\le 1, \varepsilon'>0$.

Let $G=SL(2, \C)$ considered as a real Lie group and $K=SU(2)$.
Consider the symmetric space $X=SL(2,\C)/SU(2)$. Let
$$A=\{a_t=\left (\ba{ll}e^{\frac t2}&0\\0&e^{-\frac t2}\ea\right
) |\, t\in \R\}.$$ Then $\phi_\lambda(a_t)=\frac{\sin (\lambda
t)}{\lambda\sinh t}$ and the Plancherel measure
$\mu(\lambda)=\lambda^2$ (see \cite{He}, p. 432). We define a
$K$-biinvariant function $g$ on $X$ by prescribing its spherical
Fourier transform
$\what{g}(\lambda)=\int_Gg(x)\phi_{-\lambda}(x)dx={\mathcal
F}(\psi)(\lambda)e^{-\frac{\lambda^2} 4}P(\lambda)$ for
$\lambda\in \R$ where $\psi$ is an even function in
$C_c^\infty(\R)$ with support $[-\zeta, \zeta]$ for some
$\zeta>0$, ${\mathcal F}(\psi)$ is its Euclidean Fourier transform
and $P$ is an even polynomial in $\R$. This means that $g$ is the
convolution (in $G$) of a smooth compactly supported
$K$-biinvariant function on $G$ with a (invariant) derivative of
the heat kernel of $X$. Indeed it is clear from Paley-Wiener
theorem that ${\mathcal F}(\psi)$ is also the spherical Fourier
transform of a $K$-biinvariant smooth function on $G$ supported in
a ball of radius $\zeta$.

 It follows that $g$ is a $K$-biinvariant
function of the $L^2$-Schwartz space of $G$. By the inversion
formula for the spherical Fourier transform we have,
$$\ba{lll}g(a_t)&=&C\int_{\R}\what{g}(\lambda)\phi_\lambda(a_t)d\mu(\lambda)\\
&=&\frac{C}{\sinh t}\int_\R{\mathcal
F}(\psi)(\lambda)e^{-\frac{\lambda^2} 4}\lambda
P(\lambda)\sin\lambda t\, d\lambda.\ea$$ Using Fourier inversion
on $\R$, we see that $g(a_t)=\frac{C}{\sinh t}(\psi_1*_\R h)(t)$
where $*_\R$ is the convolution in $\R$, $\psi_1$ is a derivative
of $\psi$ and hence a function in $C_c^\infty(\R)$ with support
contained in $[-\zeta, \zeta]$; $h(t)=e^{-t^2}$. An easy
computation shows that
$$|g(a_t)|\le Ce^{-t^2}e^{4\zeta t-t}\le
Ce^{-\sigma(a_t)^2}\Xi(a_t)^{1-4\zeta}.$$ If we choose $\zeta>0$
such that $l=1-4\zeta>0$, then we see that the function $g$ on $X$
satisfies:

$$|g(x)|\le Ce^{-\sigma(x)^2}\Xi^l(x)(1+\sigma(x))^M \mbox { for
all } x\in X$$ for  $M>0$ and $l\in (0,1)$ and its spherical
Fourier transform $\what{g}$ satisfies
$$|\what{g}(\lambda)|\le C'e^{-\frac {\lambda^2}
4}(1+|\lambda|)^N \mbox{ for all } \lambda\in \R$$ for some $N>0$.

Thus we can find a function $g$ on $X$  which satisfies the above
estimate for any given $l\in (0, 1)$. Suppose $\varepsilon'>0$. We
choose $l$ (that is choose $\zeta$) so that $l+\varepsilon' \ge
1$. Then it is easy to verify that $g$ satisfies the estimate
(\ref{beurling-heat-condition-sharp2}) with $c\le 1$ for any
suitable large $d$.

Now suppose $c< 1$ and $\varepsilon '\ge 0$. If $c\le 0$, the
above function $g$ clearly satisfies
(\ref{beurling-heat-condition-sharp2}). We need only therefore
consider the case when $0<c<1$. Notice that we can choose $\alpha,
\beta\in \R^+$, $\alpha<1$ and $\beta<\frac 14$ satisfying  the
constraint $4\alpha\beta=c^2$ such that the above function $g$ and
its spherical Fourier transform $\what{g}$ satisfy
$$|g(x)|\le Ce^{-\alpha\sigma(x)^2}\Xi(x) \mbox { for
all } x\in X$$  and
$$|\what{g}(\lambda)|\le C'e^{-\beta\lambda^2
}\mbox{ for all } \lambda\in \R.$$ Clearly the pair ($g,\what{g}$)
satisfy (\ref{beurling-heat-condition-sharp2}).

 From the construction of $g$ it is clear that there are infinitely many linearly independent functions
satisfying the estimate in case (ii). This example is a
modification of the example given in \cite{SiSu}.

\section{Consequences of Beurling's theorem}\newsection
In this section  we will  justify our claim made in the
introduction that this extension of the Beurling-H\"ormander
theorem is the ``master theorem'', that is all other theorems of
this genre follow from theorem 3.1. First we consider the
Gelfand-Shilov theorem.
\begin{Theorem}{\em(Gelfand-Shilov)}
\label{G-S-thm} Let $f\in  L^2(X)$. Suppose $f$ satisfies
$$\ba{ll} (1)& \int_X\frac{|f(x)|e^{\frac{(\alpha\sigma(x))^p}{p}}\Xi(x)}{(1+\sigma(x))^N}dx<\infty,\\&\\
(2)& \int_{\mathfrak a^*}\frac{\|\what{f}(\lambda)\|_2
e^{\frac{(\beta|\lambda|)^q}{q}}}{(1+|\lambda|)^N}d\mu(\lambda)<\infty,\ea$$
where $ 1 <p <{\infty}$, $\frac 1p+\frac 1q=1$ and $N$ is a
nonnegative integer.

{\em (a)} If  $\alpha\beta>1$ then $f= 0$ almost everywhere.\\

{\em (b)} If $\alpha\beta=1$ and $p\neq 2$ (and hence $q\neq 2$)
then
$f= 0$ almost everywhere.\\

{\em (c)} If $\alpha\beta=1$, $p=q=2$ and $N<d_X+1$ then $f=0$ almost everywhere\\

{\em (d)} If $\alpha\beta=1$, $p=q=2$ and $N\ge d_X+1$ then $f$ is
a $K$-finite function of the form described in theorem {\em 3.1}.
In particular if  $N= d_X+1$, then $f$ is a constant multiple of
the heat kernel $h_t$ for some $t>0$.

\end{Theorem}
\begin{proof}
(a) Since
$\frac{\alpha^p}{p}\sigma(x)^p+\frac{\beta^q}{q}|\lambda|^q\ge
\alpha\beta\sigma(x)|\lambda|$ and
$(1+\sigma(x)+|\lambda|)^{2N}\ge (1+\sigma(x))^N(1+|\lambda|)^N$,
from the assumptions (1) and (2) we obtain:
\begin{equation}\int_X\int_{\mathfrak a^*}
\frac{|f(x)\|\what{f}(\lambda)\|_2e^{\alpha\beta\sigma(x)|\lambda|}\Xi(x)}{(1+\sigma(x)+|\lambda|)^{2N}}
dxd\mu(\lambda)<\infty \label{GS-proof-1}
\end{equation}
But as $\alpha\beta>1$ we conclude that $f=0$ almost everywhere
(see section 4).

(b) Fix a $\delta\in \what{K}_0$ and an integer $j$ such that
$1\le j\le d(\delta)$. We will show that $f_{\delta, j}= 0$ almost
everywhere. Note that conditions (1) and (2) of the theorem can be
reduced (respectively) to
\begin{equation}\int_X\frac{|f_{\delta, j}(x)|e^{\frac{(\alpha\sigma(x))^p}{p}}\Xi(x)}{(1+\sigma(x))^N}dx<\infty
\label{GS-proof-2.0}
\end{equation} and
\begin{equation} \int_{\mathfrak a^*}\frac{|\what{f}_{\delta,
j}(\lambda)|
e^{\frac{(\beta|\lambda|)^q}{q}}}{(1+|\lambda|)^N}d\mu(\lambda)<\infty.
\label{GS-proof-2}
\end{equation}
Therefore we can confine ourselves to the $(\delta, j)$-th
component of the function. Using $\alpha\beta=1$  we can argue as
in (a) and show  that
\begin{equation}\int_X\int_{\mathfrak a^*}
\frac{|f_{\delta,
j}(x)|\what{f}_{\delta,j}(\lambda)|e^{\sigma(x)|\lambda|}\Xi(x)}{(1+\sigma(x)+|\lambda|)^{2N}}
dxd\mu(\lambda)<\infty \label{GS-proof-2.1}
\end{equation} and thereby conclude from theorem 3.1 that $\what{f}_{\delta,j}$ is either identically zero or
of the form $P_{\delta, j}(\lambda)e^{-\beta_0\lambda^2}$ for some
$\beta_0>0$.

Now, if  we consider the case when $1<p<2$, then we see that
unless $\what{f}_{\delta, j}=0$ almost everywhere, it cannot
satisfy (\ref{GS-proof-2}) because $q>2$.

Next we take up the case when $p>2$ and hence $1<q<2$. Since
$\mu(\lambda)$ has polynomial growth (see
(\ref{estimate-plancherel-measure})) $\what{f}_{\delta,
j}=P_{\delta, j}e^{-\beta_0\lambda^2}$ satisfies:
$$\int_{\mathfrak a^*}\frac{|\what{f}_{\delta, j}(\lambda)|
e^{\frac{(\gamma_0|\lambda|)^2}{2}}}{(1+|\lambda|)^M}d\mu(\lambda)<\infty,\mbox{
where } \gamma_0=\sqrt{2\beta_0}$$ for some suitable $M>0$. We
choose $\alpha_0$ such that $\alpha_0\gamma_0>1$. Since $p>2$ and
$f_{\delta,j}\in L^1(X)$  we see from (\ref{GS-proof-2.0}) that
$$\int_X\frac{|f_{\delta, j}(x)|e^{\frac{(\alpha_{\delta, j}\sigma(x))^2}{2}}\Xi(x)}{(1+\sigma(x))^N}dx<\infty.$$
But then from (a) it follows that $f_{\delta,j}=0$ almost
everywhere.

(c-d) By the above argument
$\what{f}_{\delta,j}(\lambda)=P_{\delta,
j}(\lambda)e^{-\beta_0\lambda^2}$. It follows from
(\ref{GS-proof-2}) with $q=2$,  that $\sqrt{2\beta_0}\ge \beta$.
But if $2\beta_0> \beta^2$ then $\alpha\sqrt{2\beta_0}>1$. On the
other hand $\what{f}_{\delta, j}$ satisfies (\ref{GS-proof-2})
with $q=2$ and with $\beta$ replaced by $\sqrt{2\beta_0}$ for a
suitably large $N$. Therefore by (a) $f_{\delta, j}=0$ almost
everywhere. Hence $\what{f}_{\delta,j}(\lambda)=P_{\delta,
j}(\lambda)e^{-\frac{\beta^2} 2\lambda^2}$. Now as noted earlier,
the Kostant polynomial $Q_\delta$ is a factor of $P_{\delta,j}$
and hence $\deg P_{\delta, j}\ge \deg Q_\delta=p_\delta$.
Therefore  only for finitely many $\delta\in \what{K}_0$,
$\what{f}_{\delta,j}$ can satisfy (\ref{GS-proof-2}). This proves
the first statement in (d). Substituting $\what{f}_{\delta, j}$
back in (\ref{GS-proof-2}) and using
(\ref{estimate-plancherel-measure}) it is now easy to verify that
if $N< 2+m_\gamma+m_{2\gamma}=d_X+1$ then $\what{f}_{\delta,
j}\equiv 0$ and if $N=2+m_\gamma+m_{2\gamma}=d_X+1$ then $\deg
P_{\delta,j}=0$ and hence
$\what{f}_{\delta,j}(\lambda)=Ce^{-\frac{\beta^2}2\lambda^2}$. But
that is possible only when $\delta$ is trivial. Indeed from
(\ref{Phi^j}) it follows that $\Phi^j_{i\rho, \delta}\equiv 0$
when $\delta\in \what{K}_0$ is nontrivial and $1\le j\le
d(\delta)$. Hence for such a $\delta$,
$\what{f}_{\delta,j}(i\rho)=0$ which is not possible if
$\what{f}_{\delta,j}(\lambda)=Ce^{-\frac{\beta^2}2\lambda^2}$.
Thus $f$ is a constant multiple of the heat kernel $h_t$ where
$t=\frac{\beta^2}2$.
\end{proof}
We will see below that the theorems of Morgan, Hardy and
Cowling-Price  follow  from the Gelfand-Shilov theorem proved
above.
\begin{Theorem}{\em (Morgan's theorem)}
 Let $f:X \rightarrow \C$ be
measurable and assume that,
$$\ba{ll} (1)& |f(x)| \leq C_1e^{-a \s(x)^p}\Xi(x)(1+\sigma(x))^n, \mbox { for all } x\in X\\&\\
(2)& \|\what{f}{(\lambda)}\|_2 \leq C_2e^{-b|\lambda|^q}, \mbox{
for all }\lambda\in \mathfrak a^*\equiv \R\ea$$ where $C_1,C_2$
and $a,b$ are positive constants, $n$ is a nonnegative integer,
$1< p < {\iy}$ and
$\f 1{p} + \f{1}{q} = 1 $.\\

{\em (a)} If $(a p )^{\f{1}{p}}(b q)^{\f{1}{q}} >1 $, then $f= 0$
almost everywhere.

{\em (b)} If $(a p )^{\f{1}{p}}(b q)^{\f{1}{q}}=1$ and $p\neq 2$
then $f= 0$ almost everywhere.

{\em (c)} If  $p=q=2$ and $(ap )^{\f{1}{p}}(b q)^{\frac{1}{q}}=1$,
that is $ab=\frac 14$  then $f$ is a constant multiple of the heat
kernel.
\end{Theorem}
\begin{proof}
Let $a=\frac{\alpha^p}p$ and $b=\frac{\beta^q}q$. Then $f$ and
$\what{f}$ satisfies theorem \ref{G-S-thm} for some suitable $N$.
The condition $(a p )^{\f{1}{p}}(b q)^{\f{1}{q}}\ge 1$ translates
as $\alpha \beta\ge 1$. Thus (a) and (b) follow from (a) and (b)
of theorem \ref{G-S-thm}. For (c) again we use the proof of (c-d)
of theorem \ref{G-S-thm}, to conclude that $\what{f}_{\delta,
j}(\lambda)=P_{\delta, j}(\lambda)e^{-b\lambda^2}$. But because of
the condition (2) of this theorem $P_{\delta, j}$ is a constant.
But this implies that only for trivial $\delta=\delta_0$,
$\what{f}_{\delta, j}$ can be nonzero and the spherical Fourier
transform of $f$ is $Ce^{-b\lambda^2}$ (see the argument at the
end of the proof of theorem \ref{G-S-thm}, (c-d)). That is $f$ is
a constant multiple of the heat kernel at $t=b$.
\end{proof}
\begin{Remark}{\em
The nonnegative integer $n$ in condition (1) of Morgan's theorem
should satisfy $n\ge \frac{m_\gamma+m_{2\gamma}-2} 2$, otherwise
the heat kernel will not be accommodated in this  inequality. That
is if we start with $n<\frac{m_\gamma+m_{2\gamma}-2} 2$, then in
case (c) also we have $f= 0$ almost everywhere. }
\end{Remark}

Morgan's theorem implies the well-known Hardy's theorem as a
particular case ($p=q=2$). To stress this point we will write it
as a separate theorem.
\begin{Theorem}{\em (Hardy's theorem)}
 Let $f:X \rightarrow \C$ be
measurable and assume that,
$$\ba{ll} (1)& |f(x)| \leq C_1e^{-a \s(x)^2}\Xi(x)(1+\sigma(x))^n,  \mbox { for all } x\in X\\&\\
(2)& \|\what{f}{(\lambda)}\|_2 \leq C_2e^{-b|\lambda|^2}, \mbox{
for all }\lambda\in \mathfrak a^*\equiv \R\ea$$ where $C_1,C_2$
and $a,b$ are positive constants, $n$ is a
nonnegative integer, $1< p < {\iy}$.\\

{\em (a)} If $ab>\frac 14$, then $f= 0$ almost everywhere.

{\em (b)} If $ab=\frac 14$ then $f$ is a constant multiple of the
heat kernel.
\end{Theorem}

\begin{Theorem} {\em (Cowling-Price)} Let $f: X \rightarrow \C$ be
measurable and assume that for positive constants $a,b$ and
nonnegative integers $m,n$,
$$\ba{ll} (1)& \int_X\frac{(|f(x)|e^{a\sigma(x)^2}\Xi(x)^{\frac
2{p_1}-1})^{p_1}}{(1+\sigma(x))^{m}}<\infty\\&\\
(2)& \int_{\mathfrak
a^*}\frac{(\|\what{f}(\lambda)\|_2e^{b|\lambda|^2})^{p_2}}{(1+|\lambda|)^n}d\mu(\lambda)<\infty,\ea$$
where $ 1 \leq p_1,p_2 <\infty$,\\

{\em (a)} If $ab>\frac 14$, then $f=0$ almost everywhere.\\

{\em (b)} If $ab=\frac 14$ then $f$ is a $K$-finite function of
the form described in theorem {\em 3.1}. In particular if
$d_X<n\le d_X+p_2$,  then $f$ is a constant multiple of the heat
kernel.
\end{Theorem}

\begin{proof}
Let us first assume $p_1$ and $p_2$ are greater than 1.
 Let $q_1$ and $q_2$ be respectively the conjugates of
$p_1$ and $p_2$, that is $\frac 1{p_i}+\frac 1{q_i}=1, i=1,2$.
Using the estimate of $\Xi(x)$ given in section 2 we note that
$\frac{\Xi(x)^{\frac 2{q_1}}}{(1+\sigma(x))^{\frac {m'}{q_1}}}$ is
in $L^{q_1}(X)$ if $m'>3$. Therefore it follows from condition (1)
in the hypothesis that
$$\ba{ll}&\int_X\frac{|f(x)|e^{a\sigma(x)^2}\Xi(x)^{\frac
2{p_1}-1}}{(1+\sigma(x))^{\frac{m}{p_1}}}\times\frac
{\Xi(x)^{\frac{2}{q_1}}}{(1+\sigma(x))^{\frac{m'}{q_1}}}dx\\&\\
=&\int_X\frac{|f(x)|e^{a\sigma(x)^2}\Xi(x)}{(1+\sigma(x))^{N_1}}dx
<\infty,\ea$$ where $N_1=\frac{m}{p_1}+\frac{m'}{q_1}$. Similarly
using (\ref{estimate-plancherel-measure}) we see that if
$n'>1+m_\gamma+m_{2\gamma}$, then
$$\int_{\mathfrak
a^*}\frac{|\what{f}(\lambda)e^{b|\lambda|^2}}{(1+|\lambda|)^{N_2}}d\mu(\lambda)<\infty$$
where $N_2=\frac{n}{p_2}+\frac{n'}{q_2}$.

When either $p=1$ or $p_2=1$ then the above two inequalities are
evident.

Thus this becomes a particular case of theorem \ref{G-S-thm}  when
$p=q=2$, $N=\max\{N_1, N_2\}$ and $a=\frac{\alpha^2}2,
b=\frac{\beta^2}2$. Note that the conditions $ab>\frac 14$ and
$ab=\frac 14$ in the hypothesis translate  as $\alpha\beta>1$ and
$\alpha\beta=1$ respectively, when we fit them in theorem
\ref{G-S-thm}. The result now follows from (a), (c) and (d) of
theorem \ref{G-S-thm} in a fashion similar  to what was used in
the previous theorems in this section. We omit the details to
avoid repetitions.
\end{proof}
In the above theorem, we may take either $p_1$ or $p_2$ or both to
be infinity. The condition (1) with $p_1=\infty$ means that
$g(x)=|f(x)|e^{a\sigma(x)^2}\Xi(x)^{-1}(1+\sigma(x)^{-m}$ is a
bounded function on $X$ for $m$ as above. Hence $\frac
{g(x)\Xi(x)^2}{(1+\sigma(x))^{N_1}}$ is integrable where
$N_1=m'>3$ as described above. That is as above
$\int_X\frac{|f(x)|e^{a\sigma(x)^2}\Xi(x)}{(1+\sigma(x))^{N_1}}
dx<\infty$.

Similarly for $p_2=\infty$  we arrive at $\int_{\mathfrak
a^*}\frac{|\what{f}(\lambda)|e^{b|\lambda|^2}}{(1+|\lambda|)^{N_2}}d\mu(\lambda)<\infty$
for $N_2=n'>1+m_\gamma+m_{2\gamma}$, since
$|\what{f}(\lambda)|e^{b|\lambda|^2}$ is bounded on $\mathfrak
a^*$. Note that the case $p_1=p_2=\infty$ of the Cowling-Price
theorem implies Hardy's theorem.

Some parts of these theorems were proved independently on
symmetric spaces. Part (a) of  Hardy's theorem was proved in
\cite{SiSu}, \cite{CSS}, \cite{EEKK}, while part (b)  was proved
in \cite{NR1}, \cite{Th1}. Part (a) of  Cowling-Price theorem  was
proved in \cite{Sengupta2} and in \cite{NR2} and  part (b) was
proved in \cite{RS}. Part (a) of  Morgan's theorem was proved in
\cite{Sengupta2}.
\section{Concluding Remarks}
Demange in his thesis (\cite{Demange}) further generalized theorem
1.2:
\begin{Theorem}{\em (Demange 2004)}
For  two nonzero functions $f_1, f_2\in L^2(\R)$, if
$$\ba{ll}(1)&\int_\R\int_\R\frac{|f_1(x)||\what{f_2}(\lambda)|e^{|x||\lambda|}}{(1+|x|+|\lambda|)^d}dxd\lambda<\infty\\&\\
(2)&\int_\R\int_\R\frac{|f_2(x)||\what{f_1}(\lambda)|e^{|x||\lambda|}}{(1+|x|+|\lambda|)^d}dxd\lambda<\infty
\ea$$ then $f_1(x)=P_1(x)e^{-\alpha x^2}$ and
$f_2(x)=P_2(x)e^{-\alpha x^2}$ for some positive constant $\alpha$
and polynomials $P_1$, $P_2$.
\end{Theorem}

A careful reader will observe that  this theorem can be extended
to symmetric spaces using our technique to have the following
interesting consequence:

We consider two rank 1 symmetric spaces $X_1=G_1/K_1$ and
$X_2=G_2/K_2$. Let $dx$ and $dy$ be  the $G_1$ and $G_2$ invariant
measures on $X_1$ and $X_2$ respectively. Let $\mu_i$ be the
corresponding  Plancherel measure for $X_i$ and let $\sigma_i$,
$\Xi_i$ be the $\sigma$ and $\Xi$ functions on $X_i$, $i=1,2$. Let
$f_1\in  L^2(X_1)$ and $f_2\in L^2(X_2)$ be two nonzero functions.
\begin{Theorem} 
Let $f_1$ and $f_2$ as above satisfy
$$\ba{ll}(1)&\int_{X_1}\int_{\mathfrak a_2^*}\frac{|f_1(x)|\|\what{f_2}(\nu)\|_2e^{\sigma_1(x)|\nu|}\Xi_1(x)}
{(1+\sigma_1(x)+|\nu|)^d}dxd\mu_2(\nu)<\infty\\ &\\
(2)&\int_{X_2}\int_{\mathfrak
a_1^*}\frac{|f_2(y)|\|\what{f_1}(\lambda)\|_2e^{\sigma_2(y)|\lambda|}\Xi_2(y)}
{(1+\sigma_2(y)+|\lambda|)^d}dyd\mu_1(\lambda)<\infty.\ea$$  Then
$f_1$ is a derivative of the heat kernel $h^1_\alpha$ of $X_1$ and
$f_2$ is a derivative of the heat kernel $h^2_\alpha$ of $X_2$ for
some instant $\alpha>0$.
\end{Theorem}

We take $X_1=X_2=X$ and obtain the following corollary.
\begin{Corollary} 
Let two nonzero functions $f_1, f_2\in L^2(X)$ satisfy
$$\ba{ll}(1)&\int_X\int_{\mathfrak a^*}\frac{|f_1(x)|\|\what{f_2}(\lambda)\|_2e^{\sigma(x)|\lambda|}\Xi(x)}
{(1+\sigma(x)+|\lambda|)^d}dxd\mu(\lambda)<\infty\\&\\
(2)&\int_X\int_{\mathfrak
a^*}\frac{|f_2(x)|\|\what{f_1}(\lambda)\|_2e^{\sigma(x)|\lambda|}\Xi(x)}
{(1+\sigma(x)+|\lambda|)^d}dxd\mu(\lambda)<\infty.\ea$$ Then $f_1$
$($respectively $f_2$$)$  is a derivative of the heat kernel
$h_\alpha$ for some instant $\alpha>0$.
\end{Corollary}

The proof of the theorem  above proceeds along entirely similar
lines to that of the proof of the main theorem of this article. We
therefore omit it.

\end{document}